\documentclass[11pt,amssymb]{amsart}
\usepackage{amssymb}
\usepackage[mathscr]{eucal}

\newcommand{\A}{{\mathbb A}}

\newcommand{\Z}{{\mathbb Z}}

\newcommand{\C}{{\mathbb C}}
\newcommand{\Q}{{\mathbb Q}}

\newcommand{\R}{{\mathbb R}}

\newcommand{\rk}{\mathrm{rk}}
\newcommand{\SL}{\mathrm{SL}}
\newcommand{\Sp}{\mathrm{Sp}}
\newcommand{\GL}{\mathrm{GL}}

\newcommand{\ct}{\mathrm{ct}}

\newcommand{\md}{{\rm mod}\ }

\newcommand{\cO}{\mathscr{O}}

\newcommand{\cG}{\mathcal{G}}

\def\ni{{\noindent}}
\begin{document}

\title[Congruence subgroup problem after
Bass-Milnor-Serre]{Developments on the congruence subgroup problem
after the work of Bass, Milnor and Serre}
\author[Prasad]{Gopal Prasad}
\author[Rapinchuk]{Andrei S. Rapinchuk}

\address{Department of Mathematics, University of Michigan, Ann
Arbor, MI 48109}

\email{gprasad@umich.edu}

\address{Department of Mathematics, University of Virginia,
Charlottesville, VA 22904}

\email{asr3x@virginia.edu}

\maketitle

\noindent {\bf 1. Introduction.} The papers of Bass-Milnor-Serre
\cite{BMS} and Serre \cite{S1}  initiated the investigation of the
congruence subgroup problem for arbitrary semi-simple algebraic
groups. The aim of this article is to survey the results in the area
obtained after, and influenced by, \cite{BMS}, and also to discuss
some remaining open questions. We would like to point out that the
impact of \cite{BMS} can be seen far beyond the congruence subgroup
problem; in particular, it and the paper \cite{St0} of Steinberg
inspired the development of algebraic $K$-theory. In this article,
however, we will limit our account exclusively to the congruence
subgroup problem for linear algebraic groups. For abelian
\mbox{varieties} the problem originated with a question of J.W.S.\,Cassels on elliptic curves,
and was settled by Serre in \cite{S2}, \cite{S3}. We refer the reader
interested in $K$-theoretic aspects to the classical monographs
\cite{B}, \cite{M}, and also to the more recent book \cite{Ros}
which includes a discussion of applications of algebraic $K$-theory
to topology.

The literature devoted to the congruence subgroup problem is quite
substantial -- the reader can find detailed references in the survey
articles \cite{Rag3}, \cite{Rag4}, \cite{Rap3}, \cite{Rap4}, and the
book \cite{Su} (the last one treats only elementary aspects of the
theory). In this article, we will focus our discussion on how the
results established in \cite{BMS} and \cite{S1} fit in, and can be
derived from, the more general results and techniques available today.

\vskip3mm

\noindent {\bf 2. The congruence subgroup problem}

\vskip3mm

In this section, we will recall the
congruence subgroup problem as originally described in \cite{BMS}, \S 14,
make some additional comments, and give the reduction of the
problem (over number fields) to absolutely simple simply connected groups.
In the next section, we will discuss the connection
between the congruence subgroup and the metaplectic problems, first pointed
out in \cite{BMS}, \S 15, in a somewhat more general framework. The
subsequent sections are devoted to the results obtained on these
problems after \cite{BMS} and \cite{S1}.
\vskip1mm

\ni{\bf 2.1.} Throughout this article, $k$ will denote a global field, i.e.,
either a number field, or the function field of an algebraic curve
over a finite field. We let $V$ denote the set of all places of $k,$
and $V_f$ (resp., $V_{\infty}$) the subset of nonarchimedean
(resp., archimedean) places, and as usual, for $v \in V$, we let
$k_v$ denote the corresponding completion. Let $S$
be a (not necessarily finite) subset of $V$ which everywhere except in  the definition of the
``$S$-metaplectic kernel" $M(S , G)$ below, will be assumed to be nonempty and
containing $V_{\infty}$ if $k$ is a number field. Then the
corresponding ring of $S$-integers is defined to be
$$
\cO_S = \{ x \in k \ \vert \ v(x) \geqslant 0 \ \ \text{for all} \ \
v \notin S \}.
$$
Now, given an algebraic $k$-group $G,$ we fix a $k$-embedding $G
\stackrel{\iota}{\hookrightarrow} {\GL}_n$ and then define the
corresponding group of $S$-integral points $\Gamma = G(\cO_S)$ to be
$G(k) \cap \GL_n(\cO_S)$ (or, more precisely,
$\iota^{-1}(\iota(G(k)) \cap \GL_n(\cO_S))$). For an ideal
$\mathfrak{a}$ of $\cO_S,$ the {\it principal $S$-congruence
subgroup} of level $\mathfrak{a}$ is defined to be
$\Gamma_{\mathfrak{a}} := \Gamma \cap \GL_n(\cO_S , \mathfrak{a})$,
where $\GL_n(\cO_S , \mathfrak{a})$ is the subgroup of $\GL_n(\cO_S)$
consisting of matrices congruent to the identity matrix modulo $\mathfrak{a}$.
A subgroup $\Gamma'$ of $\Gamma$ that contains
$\Gamma_{\mathfrak{a}}$ for some nonzero $\mathfrak{a}$ is called an
{\it $S$-congruence subgroup} (in which case obviously $[\Gamma :
\Gamma'] < \infty$). The {\it congruence subgroup problem} (CSP)
(for $\Gamma$) is the question {\it whether every subgroup of finite
index of $\Gamma$ is an $S$-congruence subgroup} (cf.\:\cite{BMS},
p.\:128). (The answer to this question is independent of the
$k$-embedding $G \stackrel{\iota}{\hookrightarrow} {\GL}_n$ chosen,
see 2.3 below.) In this (classical) formulation, CSP goes back to
the works of Fricke and Klein on automorphic functions in the 19th
century who discovered that $\Gamma = \SL_2(\Z)$ (which corresponds
to $G = {\SL}_2$ over $k = \Q$, and $S = V_{\infty}$) has plenty of
subgroups of finite index which are not congruence subgroups,
providing thereby a negative answer to the problem. The first
positive results on CSP appeared only in the early 1960s when
Bass-Lazard-Serre \cite{BLS} and Mennicke \cite{Me} solved the
problem in the affirmative for $\Gamma = \SL_n(\Z),$ $n \geqslant
3.$ The pioneering study of CSP undertaken in \cite{BMS} led to an
elegant and useful reformulation of the problem and to
the introduction of an object, called the {\it $S$-congruence
kernel}, which we will now describe. \vskip1mm

\ni{\bf 2.2.} Let $\mathfrak{N}_a$ (resp., $\mathfrak{N}_c$) be the family of all
normal subgroups of finite index of  $\Gamma$ (resp., all principal
$S$-congruence subgroups). One can introduce two topologies,
$\tau_a$ and $\tau_c,$ on $\Gamma,$ compatible
with the group structure, that admit $\mathfrak{N}_a$ and
$\mathfrak{N}_c$  as fundamental systems of neighborhoods of the
identity (these topologies are called the $S$-arithmetic and
$S$-congruence topologies, respectively), and then the affirmative
answer to the congruence subgroup problem for $\Gamma$ is equivalent
to the assertion that
\begin{equation}\tag{T}
\tau_a = \tau_c.
\end{equation}
Furthermore, one observes that $\Gamma$ admits completions
$\widehat{\Gamma}$ and $\overline{\Gamma}$ with respect to
$\tau_a$ and $\tau_c,$ and since $\tau_a$ is finer
than $\tau_c,$ there is a continuous homomorphism $\widehat{\Gamma}
\stackrel{\pi^{\Gamma}}{\longrightarrow} \overline{\Gamma}.$ Then
(T) amounts to the assertion that the {\it $S$-congruence kernel}
$$
C^S(G) := \ker \pi^{\Gamma}
$$
is the trivial group. In any case, since $\widehat{\Gamma}$ and
$\overline{\Gamma}$ can be described in terms of projective limits
of finite groups as follows
$$
\widehat{\Gamma} = \lim_{\longleftarrow}\!_{N \in \mathfrak{N}_a}\:
\Gamma/N \ \ \ \text{and} \ \ \ \overline{\Gamma} =
\lim_{\longleftarrow}\!_{\mathfrak{a} \neq 0}\:
\Gamma/\Gamma_{\mathfrak{a}},
$$
the homomorphism $\pi^{\Gamma}$ is surjective and the $S$-congruence
kernel $C^S(G)$ is a profinite group. Another useful observation is
that the families $\mathfrak{N}_a$ and $\mathfrak{N}_c$ constitute
fundamental systems of neighborhoods of the identity for topologies,
compatible with the group structure, not only on $\Gamma = G(\cO_S),$
but also on $G(k).$ We will let $\tau_a$ (resp.,\,$\tau_c$) denote
also the topology on $G(k)$ defined by the family $\mathfrak{N}_a$
(resp.,\,$\mathfrak{N}_c$). Moreover, $G(k)$ admits completions
$\widehat{G}$ and $\overline{G}$ with respect to $\tau_a$ and
$\tau_c$ respectively, and both $\widehat{G}$ and $\overline{G}$ are
topological groups; see, for example, \cite{S1}, 1.1, or the first
paragraph of \S 9 in \cite{PR1}. There is a continuous homomorphism
$\widehat{G} \stackrel{\pi}{\longrightarrow} \overline{G}.$ It
easily follows from the definitions that $\widehat{\Gamma}$ and
$\overline{\Gamma}$ can be identified with the closures of $\Gamma$
in $\widehat{G}$ and $\overline{G}$, and they are open in the
respective groups $\widehat{G}$ and $\overline{G}$, that $C^S(G) =
\ker \pi^{\Gamma}$ coincides with $\ker \pi,$ and that $\pi$ is
surjective. Thus, we arrive at the following exact sequence of
locally compact topological groups
\begin{equation}\tag{C}
1 \to C^S(G) \longrightarrow \widehat{G}
\stackrel{\pi}{\longrightarrow} \overline{G} \to 1.
\end{equation}
Notice that by construction there is a natural homomorphism $G(k)
\longrightarrow \widehat{G}$ which means that (C) splits over 
$G(k)\,(\subset \overline{G})$.  We collect
these preliminary observations in the following proposition.

\vskip2mm

\noindent {\bf Proposition 1.} {\it The $S$-congruence kernel
$C^S(G)$ is the kernel in the extension $(C)$ of locally compact
topological groups which splits over $G(k).$ It is a profinite
group which is trivial if and only if the congruence subgroup
problem for $\Gamma = G(\cO_S)$ has the affirmative answer.}

\vskip2mm

\ni{\bf 2.3.} It is not difficult to see that a $k$-isomorphism
$\varphi \colon G \to G'$ of matrix algebraic $k$-groups induces an
isomorphism $\widehat{\varphi} \colon C^S(G) \to C^S(G')$ of the
corresponding $S$-congruence kernels, for any $S,$ which in particular
implies that the answer to the congruence subgroup problem is
independent of the choice of a $k$-embedding $G
\stackrel{\iota}{\hookrightarrow} \mathrm{GL}_n.$ In general,
$C^S(G)$ provides a natural measure of deviation from the
affirmative answer to CSP for $\Gamma = G(\cO_S),$ so by introducing
the former, \cite{BMS} and \cite{S1} offered a new quantitative
formulation of the congruence subgroup problem as the problem of
computation of $C^S(G)$ for all $G$ and $S$. The results obtained in
\cite{BMS} exhibited two possibilities for $C^S(G)$ for $G =
\mathrm{SL}_n$ ($n \geqslant 3$) and $\mathrm{Sp}_{2n}$ ($n
\geqslant 2$): if $S$ contains a noncomplex place (i.e., a place $v$
such that $k_v \neq \C$), then $C^S(G)$ is trivial; otherwise
$C^S(G)$ is isomorphic to the (finite cyclic) group $\mu_k$ of all
roots of unity in $k.$ Subsequently, it was shown in \cite{S1}  that
this description of $C^S(G)$ remains valid for $G = \mathrm{SL}_2$
if $\vert S \vert
> 1;$ on the other hand, if $\vert S \vert = 1$, then $C^S(G)$ is
infinite. Since this article is a commentary on
\cite{BMS}, we will focus our attention on the results yielding the
finiteness and the precise description of $C^S(G)$ for various $G$,
and just briefly mention now some references dealing
with the structure of $C^S(G)$ when it is infinite. It was shown by
O.V.\:Melnikov that for $\Gamma = \SL_2(\Z)$ the congruence kernel is
a free profinite group of countable rank, and this result was
extended by Zalesskii \cite{Zal2} to all arithmetic lattices in
$\SL_2(\R).$ Lubotzky \cite{Lu2} showed that for $\Gamma = \SL_2(\cO)$,
where $\cO$ is the ring of integers of an imaginary quadratic number field
(such a $\Gamma$ is called a Bianchi group), the congruence kernel contains a closed
normal subgroup which is a free profinite group of countable rank.
Regarding the structure of $C^S(G)$ when it is infinite and $k$ is
of positive characteristic - see \cite{Zal1} and the recent paper 
\cite{MPSZ}.
\vskip1mm

\ni{\bf 2.4.} In \cite{Ch} Chevalley settled the congruence subgroup
problem in the affirmative for all algebraic tori over number fields
and any finite $S.$ (\cite{Ch} contains two different proofs
of this result, but it was pointed out by Serre in the 1960s that
while the first proof is perfectly correct, the second one is wrong.)

If $k$ is a number field, then the CSP for the one dimensional additive
group $\mathbb{G}_a$ is easily seen to have the affirmative
solution, for any $S.$ On the other hand, it is well-known that in
this case, any unipotent $k$-group $U$ is split (over $k$), i.e.,
there is a descending chain of normal $k$-subgroups
\begin{equation}\label{E:U}
U = U_0 \supset U_1 \supset \cdots \supset U_{r} = \{e\}
\end{equation}
such that the successive quotients $U_i/U_{i+1}$ are isomorphic to
$\mathbb{G}_a,$ for all $i = 0, \ldots , r-1.$ Then an easy
inductive argument yields the congruence subgroup property for $U,$
again for any $S.$

Now combining the affirmative solution of the CSP for tori and
connected unipotent groups, and using the Levi decomposition, one
can essentially reduce the computation of $C^S(G)$ for an arbitrary
group $G$ over a number field $k$, and any finite $S$,  to the case
where $G$ is an absolutely simple simply connected group; so we will
focus exclusively on such groups in the rest of this article. More
precisely, if $k$ is a number field, then given an arbitrary
algebraic group $G,$ one can pick a Levi $k$-subgroup $L$ of the
identity component $G^{\circ}$ and consider the derived subgroup $H
= [L , L],$ which is a semi-simple $k$-group (and in fact, it is a
maximal connected semi-simple subgroup of $G$). One proves that for
any finite $S \supset V_{\infty},$ there is a natural isomorphism
between $C^S(H)$ and $C^S(G),$ reducing the problem to $G$
semi-simple. Then one considers the universal cover $\widetilde{G}
\stackrel{\theta}{\longrightarrow} G$ defined over $k$,  and relates
$C^S(G)$ and $C^S(\widetilde{G}).$ It turns out that if the
fundamental group $F = \ker \theta$ is nontrivial, then $C^S(G)$ is
infinite for any finite $S,$ so if one is interested in the
situations where $C^S(G)$ is finite, one needs to assume that $G$ is
simply connected. Finally, as a $k$-simple simply connected group is
obtained from an absolutely simple simply connected group, defined
over a finite extension of $k$, by restriction of scalars, we are
reduced to the case of an absolutely simple simply connected $G$.
\vskip1mm

When $k$ is of positive characteristic, a full reduction of the
computation of $C^S(G)$ to the absolutely simple simply connected
case is hardly possible as even for the additive group
$\mathbb{G}_a$ the $S$-congruence kernel is infinite, and on the other
hand, general unipotent groups can be very complicated (in
particular, they may not possess a normal series like (\ref{E:U})).
However, here is one simple fact which is often used: let $G$ be a
semi-direct product over $k$ of a connected $k$-group $H$ and a
vector group $U;$ assume that $H$ acts on $U$ without any nontrivial
fixed points and that $H(\cO_S)$ is Zariski-dense in $H;$ then there
is a natural isomorphism $C^S(H) \simeq C^S(G).$ It would be
interesting to generalize this result and show, for example, that if
$U$ is a connected unipotent normal $k$-subgroup of $G$ such that
the adjoint action of $G$ on the Lie algebra of $U$ has no nonzero
fixed points, then in case $G(\cO_S)$ is Zariski-dense in $G,$ we
have an isomorphism $C^S(G/U) \simeq C^S(G).$ This would, to a large
extent, reduce the problem to the case where $G$ is pseudo-reductive (that is,
it does not contain any nontrivial connected normal unipotent
$k$-subgroups). Furthermore, it may be possible to use the recent
results on the structure of pseudo-reductive groups, obtained in
\cite{CGP}, to reduce the computation of $C^S(G)$ to reductive
groups.

\vskip4mm

\noindent {\bf 3. \parbox[t]{11cm}{\baselineskip=2mm The metaplectic
problem}}

\vskip3mm

\ni{\bf 3.1.} Henceforth, $G$ will denote an absolutely simple simply connected
$k$-group. Based on the results in \cite{BMS} and \cite{S1}, Serre
formulated (as a footnote in \cite{S1}) the following conjecture, known 
as the {\it congruence subgroup conjecture}, 
which gives a qualitative description of $C^S(G)$ for such groups.
For $S$ finite, define the $S$-rank $\mathrm{rk}_S\: G$ of $G$ as the sum
of relative ranks $\mathrm{rk}_{k_v}\,G$, of $G$ over the completions $k_v$,
for $v \in S.$ Then {\it $C^S(G)$ is finite if $\mathrm{rk}_S\: G
\geqslant 2$ and $\rk_{k_v}\,G > 0$ for all $v \in S \setminus
V_{\infty},$\footnote{The condition $\rk_{k_v}\: G > 0$ for all $v
\in S \setminus V_{\infty}$ was missing in Serre's formulation, but
one easily sees that it is necessary for $C^S(G)$ to be finite. More
precisely, if $\mathrm{rk}_S\: G > 0$ and there is $v \in S
\setminus V_{\infty}$ with $\mathrm{rk}_{k_v}\: G = 0$, then
$C^S(G)$ is infinite.} and is infinite if $\mathrm{rk}_S\: G = 1$}
(these two cases are usually referred to as the higher rank case and
the rank one case, respectively). We notice that if $\mathrm{rk}_S\:
G = 0$, then $G(\cO_S)$ is finite and therefore $C^S(G)$ is trivial,
so this case can be excluded from further consideration. We also
notice that Serre's conjecture implies that for any infinite $S$
such that $\mathrm{rk}_{k_v}\: G > 0$ for all $v \in S \setminus
V_{\infty},$ the $S$-congruence kernel $C^S(G)$ must be finite (in fact,
trivial, see the remark following Theorem 3 below).  In this section we
will discuss an approach to attack Serre's
conjecture in the higher rank case that goes back to \cite{BMS}, \S 15, and is based on
relating the $S$-congruence kernel $C^S(G)$ to the {\it
metaplectic kernel} $M(S , G),$ which we will now define (in this
definition, the set $S$ can be arbitrary, even empty).
\vskip1mm

\ni{\bf 3.2.} Let $\A_S$ denote the $k$-algebra of $S$-ad\`eles, i.e., the
restricted topological product of the completions $k_v$, for $v \notin
S$, with respect to the rings of integers $\cO_v$ of $k_v$ for nonarchimedean
$v.$ Thus,
$$
\A_S = \{ (x_v) \in \prod_{v \notin S} k_v \ \mathbf{:} \ x_v \in
\cO_v \ \text{for almost all} \ v \notin S \cup V_{\infty} \}.
$$
Then the group $G(\A_S)$ of $\A_S$-points of an algebraic $k$-group
$G$ can also be viewed as the restricted topological product of the
$G(k_v)$'s, for $v \notin S$,  with respect to the $G(\cO_v)$'s, for $v
\notin S \cup V_{\infty},$ where the groups $G(\cO_v)$ are defined
as $G(k_v)\cap \GL_n(\cO_v)$ in terms of a {\it fixed} $k$-embedding
$G \hookrightarrow \mathrm{GL}_n.$ So, whenever necessary (cf., for
example, the statement of Theorem 4), we will regard $G(k_v)$ for
any $v \notin S$ as a subgroup of $G(\A_S).$ The group $G(\A_S)$ of
$S$-ad\`eles  is a locally compact topological group for
the natural topology with a basis consisting of the sets of the form
$\Omega \times \prod_{v \notin S \cup T} G(\cO_v)$ where $T \subset
V \setminus S$ is an arbitrary finite subset containing $(V
\setminus S) \cap V_{\infty}$ and $\Omega \subset G_T := \prod_{v
\in T} G(k_v)$ is an open subset (cf.\:\cite{PlR1}, \S 5.1, for the
details). The \mbox{($S$-)metaplectic} kernel of an absolutely
simple simply connected $k$-group $G$  is defined as follows:
$$
M(S , G) = \ker\big(H^2_{\rm m}(G(\A_S)) \stackrel{\small\rm
rest}{\longrightarrow} H^2(G(k))\big),
$$
where $H^2_{\rm m}(G(\A_S))$ denotes the second cohomology group of
$G(\A_S)$ with coefficients in $I := \R/\Z$ (with the trivial action of
$G(\A_S)$ on $I$) based on {\it measurable} cochains, $H^2(G(k))$
denotes the second cohomology of $G(k)$ with coefficients in $I$ based
on abstract cochains, and the restriction map is taken relative to
the natural diagonal embedding $G(k) \hookrightarrow G(\A_S).$ It is known that for a
locally compact second countable topological group $\cG,$ the group
$H^2_{\rm m}(\cG)$ classifies topological central extensions of
$\cG$ by $I,$ and for an abstract group $\cG,$ the group $H^2(\cG)$
classifies abstract central extensions of $\cG$ by $I.$
Furthermore, it follows from the results of Wigner \cite{Wig}
that when $S \supset V_{\infty}$ (which is always the case in the
situations arising from the congruence subgroup problem), the group
$H^2_{\rm m}(G(\A_S))$ in this definition can be identified with the
cohomology group $H^2_{\rm ct}(G(\A_S))$ defined in terms of
{\it continuous} cochains with values in $I.$
\vskip1mm

\ni{\bf 3.3.} Before describing the connection between the
$S$-congruence kernel $C^S(G)$
and the metaplectic kernel $M(S , G),$ we formulate the
following important finiteness result.

\vskip2mm

\noindent {\bf Theorem 1.} {\it Let $G$ be an absolutely simple
simply connected group over a~global field $k.$ For any (possibly,
empty) set $S$ of places of $k,$ the metaplectic kernel $M(S , G)$
is finite.}

\vskip2mm

In the final form, this theorem was proved in \cite{PR1},
Theorem 2.7, as a consequence of the particular cases considered
earlier in \cite{PRag1}, Theorems~2.10 and 3.4, and \cite{Rag1.5},
Theorem 2.1, and the finiteness of $H^2_{\rm ct}(G(k_v))$ for all
nonarchimedean $v$ established in \cite{PRag2}, Theorem 9.4, if $G$
is $k_v$-isotropic and in \cite{PRag3}, Theorem 8.1, if $G$ is
$k_v$-anisotropic. A simple proof of finiteness in the case where
$k$ is a number field and $S \supset V_{\infty}$ is given
in \cite{P} using some results of M.~Lazard.

\vskip1mm

\ni{\bf 3.4.} To link $C^S(G)$ and $M(S , G),$ one observes that the
$S$-congruence topology $\tau_c$ on $G(k)$ coincides with the
topology induced on it by the diagonal embedding $G(k)
\hookrightarrow G(\A_S).$ So, since $G(\A_S)$ is locally compact,
the completion $\overline{G}$ can be identified with the closure of
$G(k)$ in $G(\A_S),$ for any $G.$ But since $G$ is assumed to be
absolutely simple and simply connected in this section, and besides
we are only interested in the situation where $\mathrm{rk}_S\: G >
0,$ or equivalently, the group $G_S = \prod_{v \in S} G(k_v)$ is
noncompact, $G$ has the {\it strong approximation property} with
respect to $S$ (cf.\:\cite{PlR1}, 7.4), i.e.,  $G(k)$ is {\it dense}
in $G(\A_S)$. Thus, in this case, $\overline{G}$ is isomorphic to
$G(\A_S),$ and the congruence subgroup sequence (C) can be rewritten
in the form
\begin{equation}\tag{C$'$}
1 \to C^S(G) \longrightarrow \widehat{G}
\stackrel{\pi}{\longrightarrow} G(\A_S) \to 1.
\end{equation}
Next, we consider the following exact sequence obtained from the Hochschild-Serre
spectral sequence for continuous cohomology with coefficients in $I:$
\begin{equation}\label{E:25}
H^1_{\rm ct}(G(\A_S)) \stackrel{\varphi}{\longrightarrow} H^1_{\rm
ct}(\widehat{G}) \longrightarrow H^1_{\rm ct}(C^S(G))^{G(\A_S)}
\stackrel{\psi}{\longrightarrow} H^2_{\rm ct}(G(\A_S)).
\end{equation}
The fact that (C$'$) splits over $G(k)$ implies that $\rm{Im}\:
\psi$ is contained in the metaplectic kernel $M(S , G).$ Moreover,
using the finiteness of $M(S , G)$ and arguing as in the proof of
Theorem 15.1 in \cite{BMS}, one shows that ${\rm Im}\: \psi = M(S ,
G).$ Thus, (\ref{E:25}) leads to the following short-exact sequence.
\begin{equation}\label{E:26}
1 \to \mathrm{Coker}\: \varphi \longrightarrow H^1_{\rm
ct}(C^S(G))^{G(\A_S)} \longrightarrow M(S , G) \to 1.
\end{equation}
It follows from the results of Margulis \cite{Mar1} that the
commutator subgroup $[G(k) , G(k)]$ is always of finite index in
$G(k),$ using which we easily show that
\begin{equation}\label{E:27}
\mathrm{Coker}\: \varphi \simeq \mathrm{Hom}(\overline{[G(k) ,
G(k)]}/[G(k) , G(k)], \ I),
\end{equation}
where $\overline{[G(k),G(k)]}$ is the closure of $[G(k),G(k)]$ in $G(k)$ in
the $S$-congruence topology. Thus, $\mathrm{Coker}\:
\varphi$ is always finite, and moreover, is trivial if $G(k)$ is
perfect, i.e., $G(k) = [G(k) , G(k)]$ (which was an assumption in
Theorem 15.1 of \cite{BMS}). In the general case, the description of
$\mathrm{Coker}\: \varphi$ is related to the following conjecture
about the normal subgroups of $G(k),$ known as the Margulis-Platonov
conjecture (MP) for $G/k$:

\vskip2mm

\begin{center}
\parbox[t]{12cm}{\it Let $G$ be an absolutely simple simply connected group
over a global field $k,$ and let $T$ be the (finite) set of all
nonarchimedean places $v$ of $k$ such that $G$ is $k_v$-anisotropic.
Then for any noncentral normal subgroup $N$ of $G(k)$, there exists
an open normal subgroup $W$ of $G_T := \prod_{v \in T} G(k_v)$ such
that $N = \delta^{-1}(W)$, where $\delta \colon G(k) \to G_T$ is the
diagonal map. In particular, if $T = \varnothing$, then $G(k)$ does
not have proper noncentral normal subgroups.}
\end{center}

\vskip2mm

\noindent If (MP) holds for $G/k,$ we will say that the normal
subgroups of $G(k)$ have the {\it standard description}. Notice that
$G(k)$ can be viewed as an $S$-arithmetic group for $S = V^K
\setminus T$ and that (MP) is equivalent to the affirmative answer
to the congruence subgroup problem in this setting. So, it is not
surprising that the assumption that (MP) holds for $G/k$ is present
in the statements,  and plays an important role in the proofs, of
practically all results on the congruence subgroup problem.
Fortunately, the truth of (MP) has already been established for a
large number of groups. We will describe the known results in this direction 
in the following two paragraphs. 
\vskip1mm

\ni{\bf 3.5.} If $G$ is $k$-isotropic, then
obviously $T = \varnothing,$ and (MP) simply asserts that $G(k)$ does not contain
any proper noncentral normal subgroups. So, in this case (MP) is
equivalent to the Kneser-Tits conjecture that $G(k) = G(k)^+$, where
$G(k)^+$ is the normal subgroup of $G(k)$ generated by the
$k$-rational points of the unipotent radicals of parabolic $k$-subgroups.
\vskip1mm

The proof of the Kneser-Tits conjecture over global fields
was recently completed by P.\:Gille, who not only settled  the case of a rank one
form of type $^2\!E_6$, which had remained open for a long time, but also
provided a uniform geometric approach for a proof of
the conjecture over general fields for a class of groups (see his
Bourbaki talk \cite{Gi}). The truth of (MP) for anisotropic groups
of type $B_n,$ $C_n,$ $D_n$ (except the triality forms
$^{3,6}\!D_4$), $E_7,$ $E_8,$ $F_4,$ and $G_2$ was established in
the 1980s (see \cite{PlR1}, Ch.\:IX). Subsequently, (MP) was also
established for all anisotropic inner forms of type $A_n.$ This
result involved the efforts of several mathematicians: first, using
the previous results of Platonov-Rapinchuk and Raghunathan
(cf.\:\cite{PlR1}, Ch.\:IX, \S 2), Potapchik and Rapinchuk \cite{PoR}
reduced (MP) to the assertion that for a finite-dimensional central
division algebra $D$ over $k,$ the multiplicative group $D^{\times}$
cannot have a quotient which is a nonabelian finite simple group,
and then Segev \cite{Seg} and Segev-Seitz \cite{SegS} proved the
truth of this assertion over arbitrary fields (see the appendix in
\cite{RS} for a detailed description of the strategy of the entire
proof). The work of Segev \cite{Seg} initiated a series of
interesting results about finite quotients of the multiplicative
group of finite-dimensional division algebras over general fields.
Of special relevance for our account are  the various versions of
the congruence subgroup theorem established in \cite{RS} and \cite{RSS} for the 
multiplicative group $D^{\times}$ of a finite dimensional division algebra $D$ 
over an arbitrary field, which led to the result that all finite quotients of
$D^{\times}$ are solvable. The techniques developed in
\cite{RS} and \cite{RSS} were used in \cite{Rap5} to give a relatively
short proof of (MP) for anisotropic inner forms of type $A_n$ in
which the use of the classification of finite simple groups is limited to the fact that
these groups are 2-generated. So, at the time of this
writing, (MP) remains open only for (most) anisotropic outer
forms of type $A_n,$ anisotropic triality forms $^{3,6}\!D_4$, and
(most) anisotropic forms of type $E_6.$

\vskip1mm

\ni{\bf 3.6.} We now observe that the description of $\mathrm{Coker}\: \varphi$
given in (\ref{E:27}) implies that it is trivial if $N = [G(k),G(k)]$ is closed in $G(k)$
in the $S$-congruence topology (this is the case if (MP) holds,
and $\mathrm{rk}_{k_v}\: G > 0$ for all $v \in
S \setminus V_{\infty},$ which is our standing assumption). Thus, the
middle term $H^1_{\rm ct}(C^S(G))^{G(\A_S)}$ of (\ref{E:26}) is
always finite, and it is  isomorphic to $M(S , G)$ if (MP) holds.
On the other hand, it is easy to
see that since $G(k)$ is dense in $\widehat{G}$,
$$
H^1_{\rm ct}(C^S(G))^{G(\A_S)} = \mathrm{Hom}_{\rm
ct}(C^S(G)/{\overline{[C^S(G) , G(k)]}} , \ I).
$$
So, in order to recover $C^S(G),$ we need to make the fundamental
assumption, which is really pivotal in the consideration of the
higher rank case of Serre's conjecture, that $C^S(G)$ is contained in the
center of $\widehat{G}$, i.e., $C^S(G)$
is {\it central}. Then
$$
H^1_{\rm ct}(C^S(G))^{G(\A_S)} = \mathrm{Hom}_{\rm ct}(C^S(G) , I),
$$
the Pontrjagin dual of (the compact abelian group) $C^S(G).$ This
discussion leads to the first assertion of the following theorem.

\vskip2mm

\noindent {\bf Theorem 2.} {\it Let $G$ be an absolutely simple
simply connected algebraic group over a global field $k,$ and let
$S$ be a nonempty set of places, containing all the archimedean ones if $k$ is a
number field, such that $\mathrm{rk}_S\: G > 0$ and
$\mathrm{rk}_{k_v}\: G > 0$ for all $v \in S \setminus V_{\infty}.$

\vskip2mm

{\rm (1)} \parbox[t]{11.5cm}{If the congruence kernel
$C^S(G)$ is central, then it is finite. In addition, if (MP) holds,
then $C^S(G)$ is isomorphic to the
Pontrjagin dual of $M(S , G).$}

\vskip1mm

 {\rm (2)} \parbox[t]{11.5cm}{Conversely, if $C^S(G)$ is
finite and (MP) holds for $G/k$, then $C^S(G)$ is central (and hence
isomorphic to the Pontrjagin dual of $M(S , G)$).}}

\vskip3mm

Assertion (2) of this theorem has the following
generalization proved in \cite{Rap2.5} and \cite{PlR2}: if (MP) holds for $G/k$ and
$\mathrm{rk}_{k_v}\: G > 0$ for all $v \in S \setminus V_{\infty}$, then the finite
generation of $C^S(G)$ as a profinite group implies that $C^S(G)$ is
central. Thus, in the presence of (MP) for $G/k,$ we have the
following interesting dichotomy: $C^S(G)$ is either central and
finite, or it is not topologically finitely generated.

\vskip1mm

We see that the finiteness of $C^S(G)$ is equivalent to its
centrality, and if it is central, then it is isomorphic to the dual of the
metaplectic kernel $M(S , G).$ Thus, the consideration of the higher
rank case of Serre's conjecture splits into two problems: proving
the centrality of the $S$-congruence kernel and computing the
metaplectic kernel ({\it ``the metaplectic problem"}). Both the problems
have generated a large number of results which we will review in
sections 4 and 5 below. It should be noted, however, that while the
metaplectic problem has been solved completely, there are important
cases in the problem of centrality (e.g., $G = \mathrm{SL}_{1,D}$,
where $D$ is a division algebra - even a quaternion division
algebra) which look inaccessible at this time.

\vskip1mm

\ni{\bf 3.7.} The sequence (C$'$) has another property which is related to the
isomorphism in Theorem 2(1). Assume that $G(k)$ is perfect, $G$ is
$k_v$-isotropic for all $v \notin S$ and the $S$-congruence kernel
$C^S(G)$ is central. Then (C$'$) has the following universal
property (cf.\:\cite{PRag1}, proof of Theorem 2.9): given a topological central
extension $$1\to {C} \longrightarrow E \longrightarrow G(\A_S) \to
1,$$ with $C$ discrete, which splits over $G(k)$, there is a unique
continuous homomorphism
$\widehat{G} \to E$ such that the following diagram commutes:
$$\begin{array}{ccccccccc} 1 & \to &  C^S(G) & \longrightarrow &
\widehat{G} &
\longrightarrow &  G(\A_S) & \to & 1 \\ & & \downarrow & & \downarrow & &\parallel & & \\
1 & \to & C & \longrightarrow & E &  \longrightarrow & G(\A_S) & \to
& 1.
\end{array}$$

\vskip1mm

\ni{\bf 3.8.} The finiteness of $C^S(G)$ has several  important consequences.
In particular, it was proved in
\cite{BMS}, Theorem 16.2, that if $k=\Q$, then the finiteness of
$C^S(G)$ implies that any homomorphism of an $S$-arithmetic subgroup
of $G(\Q)$ into $\GL_n(\Q)$ coincides, on a subgroup of finite index
of the $S$-arithmetic subgroup, with a $\Q$-rational homomorphism of
$G$ into $\GL_n.$ The results of this nature in a more general
setting were obtained in \cite{S1}, \S 2.7. These are the first
instances of ``super-rigidity" which was proved later in 1974 in a
general framework  by G.A.\:Margulis, and was used by him to derive
his celebrated arithmeticity theorem (cf.\:\cite{Mar2}).

\vskip5mm

\noindent {\bf 4. Computation of the metaplectic kernel}
 \vskip3mm

\ni{\bf 4.1.} The metaplectic kernel $M(S,G)$ was computed in \cite{BMS},
Theorem 15.3,  for $G =
\SL_n,$ $n \geqslant 3$, and $\Sp_{2n},$ $n \geqslant 2$,
using the solution of the congruence subgroup problem for
$S$-arithmetic subgroups of these groups. The subsequent work in
this direction (essentially complete by now) took a different course,
stemming from the paper of Moore \cite{Mo}, which we will sketch
below. Moore used the results of R.\:Steinberg on central extensions
of the group of rational points of a simply connected Chevalley group over an arbitrary field to link topological
central extensions of the groups of rational points of such groups
over nonarchimedean local fields with norm residue symbols of local class field theory. More
precisely, let $G$ be an absolutely simple simply connected
$k$-split group, and let $v \in V_f.$ Then to every $x_v \in
H^2_{\rm m}(G(k_v))$, there corresponds a topological central
extension
\begin{equation}\tag{$\mathcal{E}(x_v)$}
1 \to I \longrightarrow \cG_{x_v}
\stackrel{\pi_{x_v}}{\longrightarrow} G(k_v) \to 1.
\end{equation}
On the other hand, according to Steinberg \cite{St0}, the extension
$(\mathcal{E}(x_v))$ is uniquely characterized by a certain 2-cocycle
$c_{x_v} \colon k_v^{\times} \times k_v^{\times} \to I$ (called a
Steinberg 2-cocycle). Using the fact that $(\mathcal{E}(x_v))$ is a
{\it topological} central extension, Moore showed that the
associated 2-cocycle $c_{x_v}$ is of the form $c_{x_v}(* , *) =
\chi_{x_v}((* , *)_v)$ for a unique character $\chi_{x_v} \colon
\mu_{k_v} \to I$ of the group $\mu_{k_v}$ of roots of unity in
$k_v,$ where $(*,*)_v$ is the norm residue symbol on $k_v^{\times}.$
Furthermore, the correspondence $x_v \mapsto \chi_{x_v}$ defines an
injective homomorphism $\iota_v:\, H^2_{\rm m}(G(k_v)) {\longrightarrow} \widehat{\mu}_{k_v},$
where
$\widehat{\mu}_{k_v}$ is the dual  of $\mu_{k_v}.$ Moore was
able to prove that $\iota_v$ is an isomorphism if $G =
\mathrm{SL}_2;$ for arbitrary simple simply connected $k$-split groups this
was established later by Matsumoto \cite{Ma}.

Next, Moore proved that for an arbitrary finite set $S$ of places of $k$, one has the following:
\begin{equation}\tag{$\mathbb{A}$}
H^2_{\rm m}(G(\mathbb{A}_S)) = \prod_{v \notin S}
H^2_{\rm m}(G(k_v)).
\end{equation}
Now let $x \in M(S , G),$ and let
$$
1 \to I \longrightarrow \cG_x \stackrel{\pi_x}{\longrightarrow}
G(\mathbb{A}_S) \to 1
$$
be the corresponding topological central extension. According to
$(\mathbb{A}),$ we can write $x = (x_v)_{v \notin S}$, with $x_v \in
H^2_{\rm m}(G(k_v)).$ Let $\chi_{x_v} = \iota_v(x_v)$ in the above notation.
Then the Steinberg 2-cocycle corresponding to the induced extension
\begin{equation}\label{E:1}
1 \to I \longrightarrow \pi_x^{-1}(G(k))
\stackrel{\pi_x}{\longrightarrow} G(k) \to 1
\end{equation}
is
$$
c_x(\alpha , \beta) = \prod_{v \notin S} \chi_{x_v}((\alpha ,
\beta)_v) \ \ \text{for}\ \alpha , \beta \in k^{\times}.
$$
Since $x \in M(S , G),$ the extension (\ref{E:1}) splits which
imposes the following relation
\begin{equation}\label{E:2}
\prod_{v \notin S} \chi_{x_v}((\alpha , \beta)_v) = 1 \ \ \text{for
all}\ \alpha , \beta \in k^{\times}.
\end{equation}
This relation connects the norm residue symbols for different
places of $k$, and therefore can be thought of as a {\it reciprocity law.} From
global class field theory, one knows the following reciprocity law
called the {\it Artin product formula}:
\begin{equation}\tag{P}
\prod_{v\in V} (\alpha , \beta)_v^{m_v/m} = 1,
\end{equation}
where $m = \vert \mu_k \vert$ and $m_v = \vert \mu_{k_v} \vert$ (by
definition, $m_v = 0$ if $k_v = \C$). Moore proved that this
reciprocity law is actually unique in the sense that any relation
between the norm residue symbols that holds on $k^{\times} \times
k^{\times}$ is a power of $(P).$ This fact enables one to conclude
that if $S$ contains a $v$ with $k_v \neq \C,$ so that the
corresponding norm residue symbol $(* , *)_v$ is nontrivial, then in
(\ref{E:2}) all the $\chi_v$'s must be trivial. This implies that in
this case the metaplectic kernel $M(S , G)$ is trivial for any
simple simply connected $k$-split group $G.$ On the other hand, if $S$ is
totally complex, then (P) still holds if one omits the factors
corresponding to $v \in S.$ Now for a given $\chi \in
\widehat{\mu}_k$, we define $\chi_v \in \widehat{\mu}_{k_v}$ by the
formula $\chi_v(t) = \chi(t^{m_v/m})$. Assuming that the $\iota_v$'s are
isomorphisms for all $v \notin S,$ we can consider $x = (x_v) \in
H^2_{\rm m}(G(\mathbb{A}_S))$, where $x_v \in H^2_{\rm m}(G(k_v))$ is such that
$\iota(x_v) = \chi_v.$ Then (\ref{E:2}) holds, and hence, $x \in M(S , G),$
leading eventually to an isomorphism $M(S , G) \simeq
\widehat{\mu}_k.$ This completes our brief review of how the results
of Moore \cite{Mo} yield a proof of the Metaplectic Conjecture
formulated in \cite{BMS}, \S 15, for $G = \SL_2,$ and in conjunction
with the results of Matsumoto \cite{Ma}, for all simple simply connected
$k$-split groups.
\vskip2mm

\ni{\bf 4.2.} Steinberg's generators-relations approach for the study of central
extensions of split groups was generalized by Deodhar \cite{Deo} to
quasi-split groups, which enabled him to compute the metaplectic
kernel for these groups. Further results required essentially new
techniques. In \cite{PRag2}, the cohomology group $H^2_{\ct}(G(K))$ was
computed for any absolutely simple simply connected isotropic
algebraic group $G$ over a nonarchimedean local field $K$ using the
Bruhat-Tits theory and spectral sequences. (We note here that
according to a result of Wigner \cite{Wig} the natural map
$H^2_{\ct}(G(K))\rightarrow H^2_{\rm m}(G(K))$ is an isomorphism.)
It is known that any absolutely simple simply connected anisotropic group
over such a field $K$ is of
the form $G = \mathrm{SL}_{1 , D}$, where $D$ is a finite
dimensional central division $K$-algebra. For such groups,
important qualitative results about $H^2_{\ct}(G(K))$ were
established in \cite{PRag3}; recently, Ershov \cite{Ersh} has been
able to obtain an upper bound for the order of $H^2_{\ct}(G(K))$
which is sharp if $K$ is a cyclotomic extension of $\Q_p$.

Using
their local results \cite{PRag2}, Prasad and Raghunathan
\cite{PRag1} computed the metaplectic kernel for all
absolutely simple simply connected $k$-isotropic
groups. The following result for arbitrary absolutely simple simply
connected groups, which settles the metaplectic problem, was
obtained by Prasad and Rapinchuk \cite{PR1}.

\vskip2mm

\noindent {\bf Theorem 3.} {\it Let $G$ be an absolutely simple
simply connected algebraic group defined over a global field $k,$
and $S$ be a finite (possibly, empty) set of places of $k$. 
Then the metaplectic kernel $M(S , G)$ is
isomorphic to a subgroup of $\widehat{\mu}_k,$ the dual of the group
$\mu_k$ of all roots of unity in $k.$ Moreover, if $S$ contain a
place $v_0$ which is either nonarchimedean and $G$ is
$k_{v_0}$-isotropic, or is real and $G(k_{v_0})$ is not
topologically simply connected, then $M(S , G)$ is trivial.
Furthermore, $M(\varnothing, G)$ is always isomorphic to
$\widehat{\mu}_k$.}

\vskip1mm

We note that for some groups $G$ of type $A_n$ and some $S,$ the
metaplectic kernel $M(S , G)$ can be a {\it proper} subgroup of
$\widehat{\mu}_k.$ We also note that it follows from the
above theorem that $M(S , G)$ is trivial for any infinite set $S.$ Thus,
a consequence of Theorems 2 and 3 is that the truth of the higher rank case of Serre's
conjecture and that of the Margulis-Platonov conjecture imply that
$C^S(G)$ is trivial for any infinite $S.$

\vskip5mm

\noindent {\bf 5. Centrality of the $S$-congruence kernel}
 \vskip3mm

\ni{\bf 5.1.} According to Theorem 2, the higher rank case of Serre's conjecture
is equivalent to the centrality of $C^S(G)$ provided that ${\rm
rk}_S\: G > 1$, and ${\rm rk}_{k_v}\,G > 0$ for all $v \in S
\setminus V_{\infty}.$ The centrality of $C^S(G)$  has been established in a large
number of cases:
for $G = \mathrm{SL}_n$ and $\mathrm{Sp}_{2n},$ in
\cite{BMS}, \S 14, and \cite{S1}, \S 2. The result was
extended to all split (Chevalley) groups by
Matsumoto \cite{Ma}. Some nonsplit 
isotropic groups of classical types were treated by L.N.\:Vaserstein.
A proof of Serre's conjecture for general isotropic groups was given
by Raghunathan \cite{Rag1}, \cite{Rag2}. Martin Kneser was the first to prove
the centrality of $S$-congruence kernel for a $k$-anisotropic group.
He treated the spinor groups
of quadratic
forms in $n \geqslant 5$ variables in \cite{Kn}, and then Rapinchuk and Tomanov
extended his method to some other classical groups (see the
discussion below for an outline of Kneser's method and the precise
references). Rapinchuk \cite{Rap2}, \cite{Rap2.5} also considered
some exceptional groups. At the time of this writing, the centrality
of $C^S(G)$ in the higher rank case of Serre's conjecture is not
known for any anisotropic inner form, and for most of the anisotropic outer
forms, of type $A_n,$ for the anisotropic triality forms of type
$D_4,$ and for most of the anisotropic groups of type $E_6.$ One of the
important, but surprisingly difficult, open cases is where $G
= \mathrm{SL}_{1,D}$, with $D$ is a quaternion central  division
algebra over $k.$
\vskip1mm

\ni{\bf 5.2.} The results on the centrality of $C^S(G)$ mentioned above
use a variety of techniques, and we refer the reader to the surveys
\cite{Rag3}, \cite{Rag4}, \cite{Rap3} and \cite{Rap4} for the
precise formulations and some discussion of the methods involved in
various cases. Unfortunately, at this point there does not exist a
uniform strategy to attack the problem of centrality (cf., however,
the concluding paragraph of this section), so we have decided to outline
two different approaches here, and show how each of them
yields the centrality for $G = \mathrm{SL}_n,$ $n \geqslant 3,$
originally established in \cite{BMS}, Theorem 14.1. The first one, which
can be traced back to Kneser \cite{Kn}, is based on the following
two simple propositions (cf.\:\cite{Rap1}, \cite{Rap2.5}).

\vskip2mm

\noindent {\bf Proposition 2.} {\it Let $G$ be an absolutely simple
simply connected algebraic group over a global field $k$ such that
the group $G(k)$ does not contain any proper noncentral normal subgroups. Assume
that there exists a  $k$-subgroup $H$ of $G$ with the
following properties:

\vskip2mm

{\rm (1)} the natural map $C^S(H) \stackrel{\iota}{\longrightarrow}
C^S(G)$ is surjective;

\vskip1mm

{\rm (2)} there exists a nontrivial $k$-automorphism $\sigma$ of $G$
such that $\sigma \vert H = \mathrm{id}_H.$

\vskip2mm

\noindent Then $C^S(G)$ is central in $\widehat{G}.$}

\vskip3mm

\begin{proof}
Indeed, $\sigma$ induces a continuous automorphism of the topological group
$\widehat{G}$ which we will denote by $\widehat{\sigma}.$ It follows
from condition (2) that $\widehat{\sigma}$ acts trivially on the
image ${\rm Im}\: \widehat{\iota}$ of the natural map $\widehat{H}
\stackrel{\widehat{\iota}}{\longrightarrow} \widehat{G}.$ But
according to condition (1), we have the inclusion $C^S(G) \subset
{\rm Im}\: \widehat{\iota},$\, so $\widehat{\sigma}$ acts trivially on
$C^S(G).$ Then for any $g \in \widehat{G}$ and any $x \in C^S(G),$
we have
$$
x = g(g^{-1} x g) g^{-1} = g \widehat{\sigma}(g^{-1} x g) g^{-1} =
(g\widehat{\sigma}(g)^{-1}) x (g\widehat{\sigma}(g)^{-1})^{-1}.
$$
In particular, all elements of the form $g\sigma(g)^{-1}$ with $g
\in G(k),$ commute with $C^S(G)$. On the other hand, since
$\sigma$ is nontrivial, there exists $g \in G(k)$ for which
$g\sigma(g)^{-1} \notin Z(G).$ Thus, the kernel of the conjugation
action of $G(k)$ on $C^S(G)$ is a {\it noncentral} normal subgroup,
so our assumption that $G(k)$ does not contain any proper noncentral
normal subgroups implies that the action is trivial. Then the action
of $\widehat{G}$ is also trivial, hence $C^S(G)$ is central.
\end{proof}

\vskip3mm

\noindent {\bf Proposition 3.} {\it Let $G \times X \longrightarrow
X$ be a $k$-action of $G$ on an affine $k$-variety $X,$ and let $x
\in X(k).$ Assume that for every normal subgroup $N$ of $\Gamma
=G(\cO_S)$ of finite index, the orbit $N \cdot x$ is open in the
orbit $\Gamma \cdot x$ in the topology induced from
$X(\mathbb{A}_S)$ (i.e., in the $S$-congruence topology). Then for
the stabilizer $H = G_x$ of $x,$ the natural map $C^S(H) \to
C^S(G)$ is surjective.}

\vskip3mm

\begin{proof}
It is enough to show that $C^S(G)$ is contained in the closure
$\widehat{H(\cO_S)}$ of $H(\cO_S)$ in $\widehat{G}.$ Clearly,
$$
C^S(G) = \lim_{\longleftarrow}\!_{N \in \mathfrak{N}_a}\:
\overline{N}/N,
$$
where $\overline{N}$ denotes the closure of $N$ in the
$S$-congruence topology. So, to prove the inclusion $C^S(G) \subset
\widehat{H(\cO_S)},$ it suffices to prove that
$$
\overline{N} \subset H(\cO_S)N
$$
for any $N \in \mathfrak{N}_a.$ Since the action $G(\mathbb{A}_S)
\times X(\mathbb{A}_S) \longrightarrow X(\mathbb{A}_S)$ is
continuous, the openness of  $N \cdot x$ in $\Gamma \cdot x$ implies
that there exists a congruence subgroup $\Gamma_{\mathfrak{a}}$ such
that $\Gamma_{\mathfrak{a}} \cdot x \subset N \cdot x.$ Then
$\overline{N} \subset N\Gamma_{\mathfrak{a}},$ and therefore
$$
\overline{N} \cdot x \subset (N\Gamma_{\mathfrak{a}}) \cdot x = N
\cdot x,
$$
i.e., $\overline{N} \cdot x = N \cdot x.$ It follows that
$\overline{N} = (H(\cO_S) \cap \overline{N})N,$ as required.
\end{proof}

\vskip3mm

\ni{\bf 5.3.} We will now derive the centrality of $C^S(G)$ for $G =
\mathrm{SL}_n,$ $n \geqslant 3,$ using the above two propositions. We
consider the natural action of $G$ on the $n$-dimensional vector
space $V,$  and let $H$ denote the subgroup of $G$ consisting of
matrices of the form $\mathrm{diag}(A , 1)$ with $A \in
\mathrm{SL}_{n-1}.$ Then $H$ is fixed elementwise by the
automorphism $\sigma = \mathrm{Int}\: g$ of $G$, where $g =
\mathrm{diag}(1, \ldots , 1, \alpha)$, $\alpha \in k^{\times}$,
$\alpha \ne 1$. So, in view of Proposition 2, it is enough to show
that the map $C^S(H) \to C^S(G)$ is surjective. For this, we take $x
= (0, \ldots , 0, 1)$ and observe that the orbit $G(\cO_S) \cdot x$
consists of all unimodular $y = (a_1, \ldots , a_n) \in \cO_S^n$.
(Recall that $(a_1,\ldots , a_n)\in \cO_S^n$ is said to be unimodular
if the ideal generated by
$a_1,\ldots , a_n$ is $\cO_S$.)
Furthermore, given a nonzero $a \in \cO_S,$ for any unimodular
$y\in \cO_S^n$ satisfying $y \equiv x\, (\md a^2\cO_S)$, there exists $g \in
E_n(a\cO_S)$ such that $y = g x,$ where, for a nonzero ideal
$\mathfrak{a}$ of $\cO_S,$ we let $E_n(\mathfrak{a})$ denote the
normal subgroup of $G(\cO_S)$ generated by the
elementary matrices with nondiagonal entries in
$\mathfrak{a}$ (cf.\:\cite{BMS}, \S 4). Indeed, using the fact that
$a$ is invertible modulo $a_n\cO_S,$ we can find $\beta_1, \ldots ,
\beta_{n-2} \in a\cO_S$ such that $b := \beta_1 a_1 + \cdots +
\beta_{n-2} a_{n-2} + a_{n-1}$ is prime to $a_n,$ and
then there exists $g_1 \in E_n(a\cO_S)$ taking $y$ to $z = (a_1,
\ldots , a_{n-2}, b, a_n),$ and the latter vector is $\equiv x\,(\md
a^2\cO_S).$ Next, we pick $\gamma , \delta \in \cO_S$ such that
$\gamma b + \delta a_n = 1,$ and let $g_2 \in E_n(a\cO_S)$ be an
element transforming $z$ to

\vskip2mm

\noindent $(a_1 - a_1\gamma b - a_1 \delta a_n, \ldots , a_{n-2} -
a_{n-2}\gamma b - a_{n-2} \delta a_n, b, a_n) =$

\vskip1mm

\hfill $(0, \ldots , 0, b, a_n) =: w.$

\vskip0.5mm

\noindent Finally, we consider the following sequence of
transformations
$$
w \mapsto (0, \ldots , 0, a\gamma b + a \delta a_n, b, a_n)
= (0, \ldots , 0, a, b, a_n) \mapsto
$$
$$
(0, \ldots , 0, a, b - \frac{b}{a} \cdot a, a_n - \frac{a_n - 1}{a}
\cdot a) = (0, \ldots , a, 0, 1) \mapsto x.
$$
Using the fact $b \equiv 0\,(\md a^2\cO_S)$ and $a_n \equiv 1\,(\md
a^2\cO_S),$ we easily see that this sequence of transformations is
implemented by some $g_3 \in E_n(a\cO_S).$ Then $(g_3g_2g_1)y = x,$
so $g = g_1^{-1}g_2^{-1}g_3^{-1} \in E_n(a\cO_S)$ is as required.
(We notice that the element $g$ we have constructed belongs to the
subgroup (and not only to the normal subgroup) generated by the
elementary matrices with nondiagonal entries in $a\cO_S.$)

Now, if $N \subset G(\cO_S)$ is a normal subgroup of finite index, 
then there exists a nonzero ideal $\mathfrak{a}$ of $\cO_S$ such
that $N \supset E_n(\mathfrak{a})$ (cf.\:Theorem 7.5(e) in
\cite{BMS}). Pick a nonzero $a \in \mathfrak{a}.$ It follows from
the computation above that the orbit $E_n(\mathfrak{a}) \cdot x$
contains all $y \in G(\cO_S) \cdot x$ satisfying $y \equiv x\,(\md
a^2\cO_S),$ and therefore is open in $G(\cO_S) \cdot x$ in the
$S$-congruence topology. Then the orbit $N \cdot x$ is also open in
$G(\cO_S) \cdot x,$ and hence according to  Proposition 3 the map $C^S(G_x) \to
C^S(G)$ is surjective. (One can show that for any nonzero ideal
$\mathfrak{a}$ of $\cO_S,$ the orbit $E_n(\mathfrak{a}) \cdot x$ is
precisely the set of all unimodular $y$'s satisfying $y \equiv x\,(\md
\mathfrak{a})$, but the proof of this fact is
longer, see Theorem 3.3 in \cite{B}, Ch.\:V.) Furthermore, $G_x= U
\rtimes H$ where $U$ is the unipotent radical of $G_x,$ which is a
vector group. So, it follows from the discussion at the end of
section 2 that the map $C^S(H) \to C^S(G_x)$ is an isomorphism. We
obtain that $C^S(H) \to C^S(G)$ is surjective, completing the
argument.

\vskip1mm

Kneser used a similar argument in \cite{Kn} to prove the centrality of
$C^S(G)$ for $G = \mathrm{Spin}(f)$, where $f$ is a nondegenerate 
quadratic form in $n \geqslant 5$ variables.  His
argument was extended by Rapinchuk \cite{Rap1}, \cite{Rap2} (a
detailed exposition was given in \cite{Rap2.5}) and Tomanov
\cite{Tom1}, \cite{Tom2} to other classical groups defined in terms
of sesqui-linear forms of a sufficiently large dimension,
and to the groups of type $G_2.$

\vskip2mm

\ni{\bf 5.4.} Another approach to proving centrality is based on the following
theorem (announced in \cite{Rap4}) which for isotropic groups was
proved by Raghunathan in a somewhat different form and with
additional restrictions (cf. Proposition 2.14 in \cite{Rag2}). The
general case follows from our result on the congruence subgroup
property for the groups of points over semi-local rings \cite{PR1},
\S 9 (which was the first result in the investigation of the
congruence subgroup problem that does not require any case-by-case
considerations).

\vskip3mm

\noindent {\bf Theorem 4.} {\it Let $G$ be an absolutely simple
simply connected algebraic group over a global field $k$ such that
the normal subgroups of $G(k)$ have the standard description.
Assume that for every $v \notin S$, there is  a subgroup $\cG_v$ of
$\widehat{G}$ so that the following conditions are satisfied:

\vskip2mm

\ni$(i)$ \parbox[t]{12cm} {$\pi(\cG_v) = G(k_v)$ \ for all \ $v \notin S,$ where $\pi$ is as in the exact sequence $(C')$ of {3.4};}

\vskip1mm

 \ni$(ii)$ \parbox[t]{12cm}{$\cG_{v_1}$ and $\cG_{v_2}$ commute elementwise for all $v_1 ,
v_2 \notin S,$ $v_1 \neq v_2;$}

\vskip1mm

\ni$(iii)$ \parbox[t]{12cm}{the subgroup generated by the $\cG_v$, for $v \notin S,$ is
dense in $\widehat{G}.$}

\vskip2mm

\noindent Then $C^S(C)$ is central in $\widehat{G}.$}

\vskip3mm

\ni{\bf 5.5.} We will now show how Theorem 4 applies to $G = \mathrm{SL}_n,$ $n
\geqslant 3.$ For $1 \leqslant i , j \leqslant n,$ $i \neq j,$ let
$U_{ij}$ be the 1-dimensional unipotent subgroup of $G$ formed by
the elementary matrices $e_{ij}(*).$ The following commutator relation
for the elementary matrices is well-known and simple to verify:
\begin{equation}\label{E:7}
[e_{ij}(s) , e_{lm}(t)] = \left\{ \begin{array}{ll} 1, & i \neq m, j
\neq l \\ e_{im}(st), & j = l, i \neq m \\ e_{lj}(-st), & j \neq l,
i = m \end{array} \right.
\end{equation}
It is easy to see that the topologies $\tau_a$ and $\tau_c$ of
$G(k)$ induce the same topology on each $U_{ij}(k)$ (cf.\:Theorem
7.5(e) in \cite{BMS}). So, if $\widehat{U}_{ij}$ and
$\overline{U}_{ij}$ denote the closures of $U_{ij}(k)$ in
$\widehat{G}$ and $\overline{G},$ respectively, then $\widehat{G}_S
\stackrel{\pi}{\longrightarrow} \overline{G}_S$ restricts to an
isomorphism $\widehat{U}_{ij} \stackrel{\pi_{ij}}{\longrightarrow}
\overline{U}_{ij}.$ By the strong approximation property of $k^+$,
the isomorphism $e_{ij} \colon k^+ \to U_{ij}(k),$ $t \mapsto
e_{ij}(t),$ extends to an isomorphism $\mathbb{A}_S \to
\overline{U}_{ij}$ which will be denoted $\overline{e}_{ij}.$ Then
$\widehat{e}_{ij} := \pi_{ij}^{-1} \circ \overline{e}_{ij}$ is an
isomorphism between $\mathbb{A}_S$ and $\widehat{U}_{ij}$. We will 
let $\cG_v$, for $v \notin S$, denote the subgroup of $\widehat{G}$
generated by $\widehat{e}_{ij}(t)$ for all $t\in k_v \subset
\mathbb{A}_S$ and all $i \neq j.$ Clearly, the $\cG_v$'s satisfy
condition $(i)$ of Theorem 2. The closed subgroup of $\widehat{G}$
generated by the $\cG_v$, $v\notin S$, contains
$\widehat{e}_{ij}(t)$ for all $t \in k,$ hence also $G(k),$ and therefore
it coincides with $\widehat{G},$ verifying condition $(iii).$
Finally, to check $(ii),$ we observe that the density of $k$ in
$\mathbb{A}_S$ implies that (\ref{E:7}) yields a similar expression
for $[\widehat{e}_{ij}(s) , \widehat{e}_{lm}(t)]$ for any $s , t \in
\mathbb{A}_S.$ Now, if $s \in k_{v_1}$, $t \in k_{v_2}$, where $v_1
\neq v_2,$ then $st = 0$ in $\mathbb{A}_S,$ which implies that
$\widehat{e}_{ij}(s)$ and $\widehat{e}_{lm}(t)$ commute except
possibly when $l = j$ and $m = i.$ In the latter case, as $n
\geqslant 3,$ we can pick $l \neq i , j$ and then write
$\widehat{e}_{ji}(t) = [\widehat{e}_{jl}(t) ,
\widehat{e}_{li}(1_{k_{v_2}})].$ Since $\widehat{e}_{ij}(s)$ is
already known to commute with $\widehat{e}_{jl}(t)$ and
$\widehat{e}_{li}(1_{K_{v_2}}),$ it commutes with
$\widehat{e}_{ji}(t)$ as well. This shows that $\cG_{v_1}$ and
$\cG_{v_2}$ commute elementwise, which completes the argument.

\vskip2mm

\ni{\bf 5.6.} Theorem 4 enables one to establish the centrality of $C^S(G)$ in
many other situations. For example, we will now use it to prove
centrality for $G = \mathrm{SL}_2$ when $\vert S \vert > 1$, 
established first by Serre in \cite{S1}. (It is worth noting that the
argument below, unlike Serre's original argument, does not require
Tchebotarev's Density Theorem.) To keep our notations simple, we
will give the argument for $k = \Q,$ $S = \{v_{\infty} , v_p\},$ but
it extends immediately to the general situation. Let us introduce
the following notations:
$$
u^+(a) = \left(\begin{array}{cc} 1 & a \\ 0 & 1 \end{array}\right)\
, \ u^-(b) = \left(\begin{array}{cc} 1 & 0 \\ b & 1
\end{array}\right)\ , \ h(t) = \left(\begin{array}{cc}
t & 0 \\ 0 & t^{-1} \end{array}\right).
$$
We will think of $u^+,$ $u^-,$ and $h$ as parametrizations of the
subgroups $U^+,$ $U^-,$ and $H$ of upper unitriangular, lower
tringular, and diagonal matrices in $G,$ respectively. One easily
checks that for $a , b \in k$ such that $ab \neq 1$ we have the
following commutator identity
\begin{equation}\label{E:77}
[u^+(a) , u^-(b)] = u^+\left(-\frac{a^2b}{1 - ab}\right) h\left(
\frac{1}{1 - ab} \right) u^-\left(\frac{ab^2}{1 - ab} \right).
\end{equation}
As above, we let $\widehat{U}^{\pm}$ and $\overline{U}^{\pm}$ denote
the closures of $U^{\pm}(k)$ in $\widehat{G}$ and $\overline{G},$
respectively. Again, $u^{\pm}$ induce isomorphisms
$\overline{u}^{\pm} \colon \mathbb{A}_S \to \overline{U}^{\pm},$ and
$\pi \colon \widehat{G} \to \overline{G}$ restricts to isomorphisms
$\widehat{U}^{\pm} \stackrel{\pi^{\pm}}{\longrightarrow}
\overline{U}^{\pm}$ which enables us to define isomorphisms
$\widehat{u}^{\pm} := (\pi^{\pm})^{-1} \circ \overline{u}^{\pm}
\colon \mathbb{A}_S \to \widehat{U}^{\pm}.$ For $v \notin S,$  let
$\cG_{v}$ be the subgroup of $\widehat{G}$ generated by
$\widehat{u}^+(k_v)$ and $\widehat{u}^-(k_v).$ Then the subgroups
$\cG_v$ clearly satisfy conditions $(i)$ and $(iii)$ of Theorem 4,
so we only need to verify condition $(ii).$ In other words, we need
to show that for $v_1 \neq v_2,$ the subgroups
$\widehat{u}^+(k_{v_1})$ and $\widehat{u}^-(k_{v_2})$ commute
elementwise. First, we construct {\it nonzero} $a_0 \in k_{v_1}$ and
$b_0 \in k_{v_2}$ such that $\widehat{u}^+(a_0)$ and
$\widehat{u}^-(b_0)$ commute in $\widehat{G}.$ Let $q_1 , q_2$ be
the primes corresponding to $v_1 , v_2,$ and let $q_3, q_4, \ldots$
be all other primes $\neq p.$ For each $m \geqslant 2,$ we can find
$n(m)$ divisible by $m!$ so that
$$
p^{n(m)} \equiv 1\:(\md (q_1 \cdots q_m)^{2m}).
$$
Then we can write $1-p^{n(m)} = a_mb_m$ with
$$
a_m \equiv 0\:(\md (q_2 \ldots q_m)^m), \ (a_m , q_1) = 1
$$
and
$$
b_m \equiv 0\: (\md (q_1q_3 \cdots q_m)^m), \ (b_m , q_2) = 1.
$$
For a suitable subsequence $\{ m_j \},$ we have
$$
a_{m_j} \longrightarrow a_0 \in k_{v_1}^{\times} \ \ \text{and} \ \
b_{m_j} \longrightarrow b_0 \in k_{v_2}^{\times}.
$$
Using (\ref{E:77}), we obtain
$$
[u^+(a_m) , u^-(b_m)] =
u^+(-a_m^2b_mp^{-n(m)})h(p^{-n(m)})u^-(a_mb_m^2p^{-n(m)})
\rightarrow 1 \ \text{in} \ \widehat{G},
$$
because $a_m^2b_m ,\: a_mb_m^2 \longrightarrow 0$ in $\mathbb{A}_S$,
and $h(p^{n(m)}) \longrightarrow 1$ in $\widehat{G}$ as $n(m)$ is
divisible by $m!.$ Thus, $[\widehat{u}^+(a_0) , \widehat{u}^-(b_0)]
= 1.$ Now, for $t \in k^{\times},$ the automorphism $\sigma_t$ of
$G$ given by conjugation by $\mathrm{diag}(t , 1)$ extends to an
automorphism $\widehat{\sigma}_t$ of $\widehat{G}_S.$ Then
$$
1 = \sigma_t([\widehat{u}^+(a_0) , \widehat{u}^-(b_0)]) =
[\widehat{u}^+(ta_0) , \widehat{u}^-(t^{-1}b_0)]
$$
for any $t \in k^{\times},$ and since $k^{\times}$ is dense in
$k_{v_1}^{\times} \times k_{v_2}^{\times}$ by weak approximation, we
obtain that $[\widehat{u}^+(a) , \widehat{u}^-(b)] = 1$ for {\it
all}  $a \in k_{v_1},$ $b \in k_{v_2},$ as required.

\vskip2mm

\ni{\bf 5.7.} Some elaboration of the argument given in 5.5 and 5.6 allows
one to give a relatively short proof of the centrality of $C^S(G)$ for any
$k$-isotropic $G$ if $\mathrm{rk}_S\: G \geqslant 2$ (at
least when $\mathrm{char}\: k \neq 2$) - the details are given in
\cite{PR2}. Theorem 4 can also be used to simplify the original
proof of centrality for anisotropic groups of types $E_7,$ $E_8$
and $F_4$ given in \cite{Rap2}, \cite{Rap2.5}. 

Finally, we would
like to mention some results on the centrality of $C^S(G)$ for
infinite $S$ (and then in these situations $C^S(G)$ is actually
trivial - cf.\:the remark following Theorem 3). Assume that (MP) holds for $G/k$
and that ${\rm rk}_{k_v}\: G > 0$ for all $v \in S
\setminus V_{\infty}.$ Then, as we have already mentioned prior to
the statement of Theorem~4, $C^S(G)$ is known to be central (hence
trivial) in the semi-local case, i.e., when $V \setminus S$ is
finite. A much stronger fact is established in \cite{PR2}: $C^S(G)$
is central whenever $S$ almost contains a ``generalized arithmetic
progression" (with some minor restrictions in the case where $G$ is
an outer form over $k$). These results are proved using techniques
that do not require any case-by-case considerations, hence may lead
to a general approach to Serre's congruence subgroup conjecture.

\vskip3mm

\noindent {\bf 6. Bounded generation of arithmetic subgroups and the CSP}

\vskip3mm

\ni{\bf 6.1.} Let $k$ be a number field and $\cO$ be the ring of
integers of $k$. Using some techniques developed in \cite{BMS}
(Mennicke symbols), Carter and Keller \cite{CK} established the
following remarkable fact: {\it for $n \geqslant 3,$ there exists $d
\in \mathbb{N}$ such that every $x \in \SL_n(\cO)$ is a product of
$\leqslant d$ elementary matrices} (notice that Corollary 4.3(a) in
\cite{BMS} states that every $x \in \SL_n(\cO)$ is a product of some
{\it unspecified} number of elementary matrices). This result
motivated the following.

\vskip2mm

\noindent {\bf Definition.} An abstract (discrete) group $\Gamma$ is
{\it boundedly generated} if
there exist $\gamma_1, \ldots , \gamma_t \in \Gamma$ such that
$\Gamma = \langle \gamma_1 \rangle \cdots \langle \gamma_t \rangle$,
where $\langle \gamma_i \rangle$ is the cyclic subgroup generated by
$\gamma_i.$
\vskip1mm

A profinite group $\Delta$ is boundedly generated (as a profinite
group) if there are $\delta_1, \ldots , \delta_t \in \Delta$ for
which $\Delta = \overline{\langle \delta_1 \rangle} \cdots
\overline{\langle \delta_t \rangle}$, where $\overline{\langle
\delta_i \rangle}$ is the closure of $\langle \delta_i \rangle.$

\vskip2mm

The result of Carter-Keller implies that $\Gamma = \SL_n(\cO)$ for
$n \geqslant 3$ is boundedly generated (of course, $\SL_2(\Z),$
being virtually free, is not boundedly generated). More importantly,
the use in \cite{CK} of techniques developed to solve the congruence
subgroup problem was far from coincidental: as the following theorem
shows, bounded generation of $S$-arithmetic groups always implies
that $C^S(G)$ is finite.

\vskip2mm

\noindent {\bf Theorem 5.} (\cite{Lu3}, \cite{PlR2}). {\it Let $G$ be
an absolutely simple simply connected algebraic group over a number
field $k.$ Set $\Gamma = G(\cO_S),$ and assume that the normal
subgroups of $G(k)$ have the standard description. Then bounded
generation of the profinite completion $\widehat{\Gamma}$ is
equivalent to the centrality (hence finiteness) of the
$S$-congruence kernel $C^S(G).$ In particular, if \,$\Gamma$ is
boundedly generated as a discrete group, then $C^S(G)$ is finite.}

\vskip2mm

This theorem in conjunction with \cite{CK} provides yet another way
to prove the centrality of the congruence kernel in \cite{BMS}, and
the hope is that this approach may lead to the resolution of some
new cases in the congruence subgroup problem. At the time of this
writing (2008), bounded generation of $\Gamma = G(\cO_S)$ is known
if $G$ is either split or quasi-split over $k$, and has $k$-rank $>
1$ (Tavgen \cite{Tav}),  $G = {\SL}_2$ and $\cO_S$ has infinitely
many units or equivalently ${\rm rk}_S\: G > 1$ - see \cite{Mor},
and $G = \mathrm{SO}(f)$, where $f$ is a nondegenerate quadratic form in $n
\geqslant 5$ variables over $k,$ and either the Witt index of $f$ is
$\geqslant 2,$ or it is 1 and $S$ contains a nonarchimedean place
\cite{ER}. Unfortunately, we still do not have a single example of a
$k$-anisotropic group with an infinite boundedly generated
$S$-arithmetic subgroup, nor do we know that $\SL_n(\cO_S),$ where
$n \geqslant 3,$  can be boundedly generated by a (finite) system of
generic semi-simple elements (here we call an element $x \in
\SL_n(k)$ generic if the Galois group over $k$ of its characteristic
polynomial is the symmetric group $S_n$). In view of the
significant difficulties associated
with verifying bounded generation, some efforts were made to find
other conditions that still imply the congruence subgroup property
but are easier to check. For example, in \cite{PlR2} such conditions
were formulated in terms of polynomial index growth and its
variants. Furthermore, it was shown that one of those variants can
be verified for $\SL_n,$ $n \geqslant 3,$ by a simple computation
based exclusively on the commutation relations (\ref{E:7}) for
elementary matrices. Thus, the centrality in \cite{BMS} can be
derived by methods of combinatorial group theory (cf.\:also
Steinberg \cite{St}).

\vskip1mm

\ni{\bf 6.2.} Groups with bounded generation property 
arise also in other areas. For example, in the theory of pro-$p$
groups, the bounded generation property characterizes
analytic pro-$p$ groups \cite{DdSMS}. In representation theory, one
proves that discrete groups $\Gamma$ with bounded generation
satisfying condition (Fab) (which means that any subgroup $\Delta
\subset \Gamma$ of finite index has finite abelianization
$\Delta^{\small \mathrm{ab}} = \Delta/[\Delta , \Delta]$), are
representation rigid, i.e., have only finitely many inequivalent
irreducible representations in each dimension. Bounded
generation of $\SL_n(\cO)$ with respect to elementary matrices was used by
Shalom \cite{Shal} to estimate the Kazhdan constant for this group -
notice that Shalom's methods apply to general rings, however bounded
generation has not been established yet over any ring other than
the rings of algebraic integers. In fact, it is known that
$\SL_3(\C[x])$ does not have bounded generation with respect to
elementary matrices  (\cite{vdK}). Whether or not $\SL_n(\Z[x])$ and $\SL_n(\Q[x])$, $n
\geqslant 3,$ are boundedly generated remains an open question.

\vskip5mm

\noindent {\it Acknowledgements.} It is a pleasure to thank Brian Conrad and Jean-Pierre Serre for their comments and corrections.

Both the authors acknowledge
partial support from the NSF (grants DMS-0653512 and
DMS-0502120), BSF (grant 2004083), and the Humboldt Foundation.
During the preparation of this article they enjoyed the hospitality
of the SFB 701 (Universit\"at Bielefeld).

\bibliographystyle{amsplain}

\end{document}